\documentclass[a4paper,12pt,intlimits,oneside]{amsart}

\usepackage{enumerate}

\usepackage{latexsym,amssymb}

\newcommand{\comment}[1]{}

\newcommand{\bC}{{\mathbb C}}

\newcommand{\bR}{{\mathbb R}}

\newcommand{\bZ}{{\mathbb Z}}

\def\C{{\mathcal C}}
\def\H{{\mathcal H}}

\def\P{{\mathcal P}}

\def\M{{\mathcal M}}

\def\bmo{{\mathfrak {bmo}}}
\def\weak{{\mbox{\small\rm  weak}}}
\def\h{\mathfrak h}
\def\lip{\Lambda_{\gamma}(\bR^{n})}

\def\la{\langle}
\def\ra{\rangle}

\def\BMO{\rm{B\! M\! O}}

\def\lmo{{\mathfrak {lmo}}}

\newcounter{rea}
\newcounter{rej}
\newcounter{res}
\newtheorem{thm}{Theorem}[section]
\newtheorem{prop}[thm]{Proposition}

\newtheorem{lem}[thm]{Lemma}
\newtheorem{defn}[thm]{Definition}
\newtheorem{remark}[thm]{Remark}
\begin{document}

\title[]{Products of functions in Hardy and Lipschitz or $\BMO$ spaces}
\author[A. Bonami]{Aline Bonami}
\address{MAPMO-UMR 6628,
D\'epartement de Math\'ematiques, Universit\'e d'Orleans,
45067 Orl\'eans Cedex 2, France}
\email{{\tt Aline.Bonami@univ-orleans.fr}}
\author[J. Feuto]{Justin Feuto}
\address{Laboratoire de Math\'ematiques Fondamentales, UFR Math\'ematiques et Informatique, Universit\'e de Cocody, 22 B.P 1194 Abidjan 22. C\^ote d'Ivoire}
\email{{\tt justfeuto@yahoo.fr}}

\subjclass{42B30, 42B35} \keywords{Orlicz spaces, atomic
decomposition, Hardy spaces, local Hardy spaces, bounded mean
oscillation.}
\thanks{Part of this work was done while the second author was visiting MAPMO at Orl\'eans, with a financial support of AUF (Agence Universitaire de la Francophonie).}
\begin{abstract}
We define as a distribution the product of a function (or
distribution) $h$ in some Hardy space $\H^p$ with a function $b$
in the dual space of $\H^p$. Moreover, we prove that the product
$b\times h$ may be written as the sum of an integrable function
with a distribution that belongs to some Hardy-Orlicz space, or to
the same Hardy space $\H^p$, depending on the values of $p$.
\end{abstract}
\maketitle

\section{Introduction}
For $p$ and $p'$ two conjugate exponents, with $1<p<\infty$, when
we consider two functions $f\in L^p(\bR^n)$ and $g\in
L^{p'}(\bR^n)= \left(L^p(\bR^n)\right)^*$, their product $fg$ is
integrable, which means in particular that their pointwise product
gives rise to a distribution. When $p=1$, the right substitute to
Lebesgue spaces is, for many problems, the Hardy space
$\mathcal{H}^1(\bR^n)$, whose dual is the space $\BMO(\bR^n)$. So
one may ask what is the right definition of the product of $h\in
\mathcal{H}^1(\bR^n)$ and $b\in \BMO(\bR^n)$. In this context, the
pointwise product  is not integrable in general. In order to get a
distribution, one has to define the product in a different way.
This question has been considered by the first author in a joint
work with T. Iwaniec, P. Jones and M. Zinsmeister in \cite{BIJZ}.
The present paper explores the same problem in different spaces.

The duality bracket $\langle b, h\rangle$ may be written through
the almost everywhere approximation of the factor $ b\in
\BMO(\bR^n)\,$,
\begin{equation}\label{approx}
\langle b, h\rangle =\lim_{k\to\infty}\int_{\bR^n} b_k(x)
h(x)\,dx\;,
\end{equation}
where $b_k$ is a sequence of bounded functions, which is bounded
in the space $\,\BMO(\mathbb R^n)$ and converges to $b$ almost
everywhere. For example, we can choose
\begin{equation}
b_k(x)= \left\{\begin{array}{lll} \;\;\;k & \mbox{if} &
\;\;\;k\leq  b(x)\\
b(x) & \mbox{if} & -k\leq  b(x)\leq \;\;k\\
-k & \mbox{if} & \quad\;\;\;\;\;\;  b(x)\leq -k
\end{array}\right. .\label{approx-1}
\end{equation}

We then define the product $b\times h$ as the distribution whose
action on the test function $\varphi$ in the Schwartz class, that
is $\varphi\in \mathcal{S}(\mathbb R^n)$, is given by
\begin{equation}\label{dual}
\langle b\times h, \varphi\rangle:= \langle b\varphi, h\rangle.
\end{equation}
 We use the fact that the multiplication by
$\varphi$ is a bounded operator on $\BMO(\bR^n)$. So the right
hand side makes sense in view of the duality $\H^1$-$\BMO$.
Alternatively, the Schwartz class is contained in the space of
multipliers of $\BMO(\bR^n)$, which have been studied and
characterized, see \cite{S} and the discussion below. It follows
from \eqref{approx}, used with the sequence $b_k$ given in
\eqref{approx-1}, that the distribution $b\times h$ is given by
the function $bh$ whenever this last one is integrable.

A more precise description of products $b\times h$ has been given
in \cite{BIJZ}. Namely, all such distributions are sums of a
function in $L^1(\bR^n)$ and a distribution in a Hardy-Orlicz
space $\mathcal{H}^\Phi_w$, where $w$ is a weight which allows a
smaller decay at infinity and $\Phi$ is given below. We will
consider a slightly different situation by replacing the space
$\BMO(\bR^n)$ by the smaller space $\bmo(\bR^n)$, defined as the
space of  locally integrable functions $b$ such that
\begin{equation}\label{bmo}
\sup_{|B|\leq 1}\left( \frac 1 {|B|}\int_B
|b-b_B|dx\right)<\infty\ \ \ \ \mbox{\rm and } \sup_{|B|\geq
1}\left( \frac 1 {|B|}\int_B |b|dx\right)<\infty.
\end{equation}
Here $B$ varies among all balls of $\bR^n$ and $|B|$ denotes the
measure of the ball $B$. Also $b_B$ is the mean of $b$ on the ball
$B$. Recall that the $\BMO$ condition reduces to the first one,
but for all balls and not only for balls $B$ such that $|B|<1$. We
clearly have $\bmo\subset \BMO$.
 We have the
following, which is new compared to \cite{BIJZ}.
\begin{thm}\label{main0}
For $h$ a function in $\mathcal H^1(\bR^n)$ and $b$ a function in
$\bmo (\bR^n)$, the product $b\times h$ can be given a meaning in
the sense of distributions. Moreover, we have the inclusion
\begin{equation}\label{inclusion0}
  b\times h\in L^1(\bR^n)+ \mathcal H_*^\Phi(\bR^n).
\end{equation}
\end{thm}
$\mathcal H_*^\Phi(\bR^n)$ is a variant of  the Hardy-Orlicz space related to the
Orlicz function
\begin{equation}\label{orl-log}
\Phi(t):=\frac t{\log (e+t)},
\end{equation}
 which is defined in Section 3. It contains $\mathcal H^\Phi(\bR^n)$ and is contained in the weighted Hardy-Orlicz space that has been considered in \cite{BIJZ}  for the general case of
$f\in\BMO$.

 The aim of this paper is to give some extensions of the previous
situation. Indeed, \eqref{dual} makes sense in other cases. First,
the space $\bmo$ is the dual of the local Hardy space, as proved
by  Goldberg  \cite{G} who introduced it. So it is natural to
extend the previous theorem to functions $h$ in this space, which
we do. Next, we can consider the Hardy space
$\mathcal{H}^p(\bR^n)$, for $p<1$, and its dual the homogeneous
Lipschitz space $\dot{\Lambda}_\gamma (\bR^n)$, with $\gamma:=
n(\frac 1p-1) $. Indeed, a function in the Schwartz class is also
a multiplier of the Lipschitz spaces. Our statement is
particularly simple when $b$ belongs to the non homogeneous
Lipschitz space ${\Lambda}_\gamma (\bR^n)$.

\begin{thm}\label{main}
Let $p<1$ and  $\gamma:=n(\frac 1p-1)$. Then, for $h$ a function
in $\mathcal H^p(\bR^n)$ and $b$ a function in $\Lambda_\gamma
(\bR^n)$, the product $b\times h$ can be given a meaning in the
sense of distributions. Moreover, we have the inclusion
\begin{equation}\label{inclusion}
  b\times h\in L^1(\bR^n)+ \mathcal H^p(\bR^n).
\end{equation}
\end{thm}

Again the space ${\Lambda}_\gamma (\bR^n)$ is the dual of the
local version of the Hardy space $\mathcal H^p(\bR^n)$. We will
adapt the theorem to $h$ in this space.

Let us explain the presence of two terms in the two previous
theorems. The product looses the cancellation properties of the
Hardy space, which explains the term in $L^1$. Once we have
subtracted some function in $L^1$, we recover a distribution of a
Hardy space. For $p=1$, there is a loss, due to the fact that a
function in $\bmo$ is not bounded, but uniformly in the
exponential class on each ball of measure $1$. This explains that
we do not find a function in $\H^1$, but in the Hardy-Orlicz
space.

As we will see, the proof uses a method that is linear in $b$, not
in $h$.  As in the case $\mathcal{H}^1$-$\BMO$ (see \cite{BIJZ}),
one would like to know whether the decomposition of $b\times h$ as
a sum of two terms can be obtained through linear operators, but
we are very far from being able to answer this question.

\medskip

All this study is reminiscent of problems related to commutators
with singular integrals, or Hankel operators. In particular, such
products arise when developing commutators between the
multiplication by $b$ and the Hilbert transform and looking
separately at each term. It is well known that the commutator
$[b,H]$ maps $\mathcal{H}^p(\bR)$ into
$\mathcal{H}^1_{\weak}(\bR)$ for $b$ in the Lipschitz space
$\Lambda_\gamma (\bR)$ (see \cite{J}), which means that there are
some cancellations between terms, compared to our statement which
is the best possible for each term separately. One can also
consider products of holomorphic functions in the corresponding
spaces when $\bR^n$ is replaced by the torus, considered as the
boundary of the unit disc. Statements and proofs are much simpler
and there are converse statements, see \cite{BIJZ} for the case
$p=1$,  and also to \cite{BG} where the problem is treated in
general for holomorphic functions in Hardy-Orlicz spaces in convex
domains of finite type in $\bC^n$. These results allow to
characterize the classes of symbols for which Hankel operators are
bounded from some Hardy-Orlicz space larger than $ \mathcal{H}^1$
into $ \mathcal{H}^1$.

\medskip

Another possible generalization deals with spaces of homogeneous
type instead of $\bR^n$. Since the seminal work of Coifman and
Weiss \cite {CW1, CW2}, it has been a paradigm in harmonic
analysis that this is the right setting for developing
Calderon-Zygmund Theory. The contribution of Carlos Segovia,
mainly in collaboration with Roberto Mac\'ias, has been
fundamental to develop a general theory of Hardy and Lipschitz
spaces. We will rely  on their work in the  last section, when
explaining how  properties of products of functions in Hardy and
Lipschitz or $\BMO$ spaces can generalize in this general setting.
Remark that the boundary of pseudo-convex domain in
$\mathbb{C}^n$, with the metric that is adapted to the complex
geometry (see for instance \cite{McN}), gives a fundamental
example of such a space of homogeneous type. Calder\'on-Zygmund
theory has been developed in this context, see \cite{KL} for
instance, in relation with the properties of holomorphic
functions, reproducing formulas, Bergman and Szeg\"o projections.
Many recent contributions have been done in $\H^p$ theory on
spaces of homogeneous type. We refer to \cite{GLY} and the
references given there. Tools developed by Mac\'ias, Segovia and
their collaborators play a fundamental role, like, for instance,
for the atomic decomposition of Hardy-Orlicz spaces given by
Viviani (see \cite{V} and \cite{BG}).

\medskip

\noindent{\bf Aknowledgement.} This paper is dedicated to the
memory of Carlos Segovia.

\section{Prerequisites on Hardy  and Lipschitz spaces}

 We recall here the definitions and properties that we will use later on.We follow the book of Stein
 \cite{St}.

 Let us first recall the definition of the  maximal  operator used for the definition of Hardy
 spaces.
 We fix a function $\varphi\in\mathcal S (\bR^{n})$ having integral $1$ and support in $\{|x|<1\}$. For
  $f\in\mathcal S' (\bR^{n})$ and $x$ in $\bR^{n}$, we put
\begin{equation}
 f\ast\varphi(x):=\left\langle f,\varphi(x-\cdot)\right\rangle,
\end{equation}
 and define the maximal function $\M_{\varphi}f$ by
\begin{equation}\label{maximal}
\M_{\varphi}f(x):=\sup_{t>0}\left|\left(f\ast\varphi_{t}\right)(x)\right|,
\end{equation}
where $\varphi_{t}(x)=t^{-n}\varphi\left(t^{-1}x\right)$. We also
define the truncated version of the maximal function, namely
\begin{equation}\label{maximal-tr}
\M_{\varphi}^{(1)}f(x):=\sup_{0<t<1}\left|\left(f\ast\varphi_{t}\right)(x)\right|.
\end{equation}
For $p>0$,  a tempered distribution $f$ is said to belong to the
Hardy space $\H^{p}(\bR^{n})$ if
\begin{equation}\label{def}
\|f\|_{\H^{p}(\bR^{n})}:=\left
(\int_{\bR^n}\M_{\varphi}f(x)^{p}dx\right)^{\frac 1p}<\infty.
\end{equation}
 The localized versions of Hardy spaces are
defined in the same spirit, with the truncated maximal function in
place of the maximal function. Namely, a tempered distribution $f$
is said to belong to the space $\h^{p}(\bR^{n})$ if
\begin{equation}\label{def-h}
\|f\|_{\h^{p}(\bR^{n})}:=\left
(\int_{\bR^n}\M_{\varphi}^{(1)}f(x)^{p}dx\right)^{\frac 1p}<\infty.
\end{equation}
Recall that, up to equivalence of corresponding norms, the space
$\H^p(\bR^n)$  (resp. ${\h^{p}(\bR^{n})}$) does not depend on the
choice of the function $\varphi$. So, in the sequel, we shall use
the notation $\M f$ instead of $\M_{\varphi}f$ (resp. $\M^{(1)}f$
instead of $\M_{\varphi}^{(1)}f$.

Hardy-Orlicz spaces are defined in a similar way. Given a
continuous function ${\mathcal P} :[0,\infty)\to [0, \infty)\,$
increasing from zero to infinity (but not necessarily convex,
$\mathcal P$ is called the Orlicz function), the Orlicz space $ L^
{\mathcal P}$ consists of measurable functions $f$ such that
\begin{equation}\label{orlicz}
\|f\|_{L^\P}:=\inf \left\{k >0\;; \int_{\bR^n}{\mathcal P}
 \left( k^{-1} |\,f|\,\right)\,d x \leqslant 1\right\}
<\infty\,.
\end{equation}
Then $\H^{\mathcal P}$ (resp. $\h^{\mathcal P}$) is the space of
tempered distributions $f$ such that $\M f$ is in $L^{\mathcal P}$
(resp. $\M^{(1)} f$ is in $L^{\mathcal P}$). We will be
particularly interested by the choice of the function $\Phi$ given
in \eqref{orl-log} as the Orlicz function. It is easily seen that
the function $\Phi$ is equivalent to a concave function (take
$t/(\log (c+t))$, for $c$ large enough). So there is no norm on
the space $ L^ {\Phi}$. In general, $\|\cdot\|_{L^\P}$ is
homogeneous, but is not sub-additive. Nevertheless (see
\cite{BIJZ}),
\begin{equation}\label{add}
  \|f+g\|_{L^\Phi}\leq 4 \left(\|f\|_{L^\Phi}+\|g\|_{L^\Phi}\right).
\end{equation}

\begin{defn} $L^\Phi_*$ is the space of functions $f$ such that
$$\|f\|_{L^\Phi_*}:=\sum_{j\in \bZ^n}\|f\|_{L^\Phi(j+\mathbb Q)}<\infty,$$
where $\mathbb Q$ is the unit cube centered at $0$.
\end{defn}

 We accordingly define
$\H^{\Phi}_*$ (resp. $\h^{\Phi}_*$). Using the concavity described above, we have $\Phi(st)\leq C s\Phi(t)$ for $s>1$. It follows that $L^\Phi$ is contained in $L^\Phi_*$ as a consequence of the fact that $\|f\|_{L^\Phi (j+\mathbb Q)}\leq \int_{j+\mathbb Q} \Phi(|f|)dx$. The converse inclusion is not true.

\medskip


 We
will restrict to $p\leq 1$, since otherwise Hardy spaces are just
Lebesgue spaces. We will  need the atomic decompositions of the
spaces $\H^{p}({\bR^n})$ (resp. $\h^{p}({\bR^n})$), which we
recall now.
\begin{defn}\label{atome}
Let $0<p\leq 1 < q\leq\infty,\;p<q$, and $s$ an integer. A
$(p,q,s)$-atom related to the ball $B$ is a function $a\in
L^{q}(\bR^{n})$ which satisfies the following conditions:
\begin{equation}\label{support}
  \mbox{\rm support}(a)\subset B \ \ \ \ \mbox{\rm and} \ \ \ \|a\|_q\leq |B|^{\frac
  1q-\frac{1}{p}},
\end{equation}
\begin{equation}\label{moment}
\int_{\bR^{n}}a(x)x^{\alpha}dx=0 \ \ \ ,\ \ \ \ \mbox{\rm for}\ \
0\leq\left|\alpha\right|\leq s.
\end{equation}
\end{defn}
Here $\alpha$ varies among multi-indices,  $x^\alpha$ denotes the
product $x_1^{\alpha_{1}}\ldots x_n^{\alpha_{n}}$ and
$|\alpha|:=\alpha_1+\cdots+\alpha_n$. Condition \eqref{moment} is
called the moment condition.

The atomic decomposition of $\H^p(\bR^n)$ is as follows. Let us
fix $q>p$ and $s>n\left(\frac{1}{p}-1\right)$. Then a tempered
distribution $f$ is in $\H^p(\bR^n)$ if and only if there exists a
sequence of $(p,q,s)$-atoms $a_j$ and constants $\lambda_j$ such
that
\begin{equation}\label{at-dec}
    f:=\sum^{\infty}_{j=1}\lambda_{j}a_{j}\ \ \ \mbox{ and} \ \ \
    \sum^{\infty}_{j=1}\left|\lambda_{j}\right|^{p}<\infty,
    \end{equation}
where the first sum is assumed to converge in the sense of
distributions. Moreover, $f$ is the limit of partial sums in
$\H^p(\bR ^n)$, and $\|f\|_{\H^p(\bR^n)}$ is equivalent to the
infimum, taken on all such decompositions of $f$, of the
quantities
$\left(\sum^{\infty}_{j=1}\left|\lambda_{j}\right|^{p}\right)^{\frac
1p} $.

 For the local version, we consider other kinds of atoms when the
 balls $B$ are large. We have the following, where we have fixed
 $q>p$ and $s>n\left(\frac{1}{p}-1\right)$. A tempered
distribution $f$ is in $\h^p(\bR^n)$ if and only if there exists a
sequence of functions $a_j$, constants $\lambda_j$ and balls $B_j$
for which \eqref{at-dec} holds, and such that
\begin{enumerate}
    \item [(i)] when $|B_j|\leq 1$, then $a_j$ is a $(p,q,s)$-atom related to
    $B_j$;
 \item[(ii)] when $|B_j|> 1$, then $a_j$ is supported in $B_j$ and
 $$\|a_j\|_q\leq |B_j|^{\frac 1q-\frac{1}{p}}.$$
\end{enumerate}
In other words, one still has the atomic decomposition, except
that for large balls one does not ask for any moment condition on
atoms.

\medskip

Next, let us define Lipschitz spaces. For $\delta\in \bR^{n}$  we
note $D^{1}_{\delta}=D_{\delta}$ the difference operator, defined
by setting $ D_{\delta}f(x)=f(x+\delta)-f(x)$ for $f$ a continuous
function (see \cite{g} for instance). Then, by induction, we
define $ D^{k+1}_{h}f=D_{\delta}\left(D^{k}_{\delta}f\right) $ for
$k$ a non negative integer,  so that
\begin{equation}
D^{k}_{\delta}f(x)=\sum^{k}_{s=0}(-1)^{k+s}\left(
\begin{array}{c}
    k\\ s
\end{array}\right)f(x+s\delta).
\end{equation}
For $\gamma >0$ and  $k=\lfloor \gamma\rfloor$ the integer part of
$\gamma$, we set
\begin{equation}
\left\|f\right\|_{\Lambda_{\gamma}}=\left\|f\right\|_{L^{\infty}}
+\sup_{x\in\bR^{n}}\sup_{h\in\bR^{n}\setminus\left\{0\right\}}
\frac{\left|D^{k+1}_{\delta}f(x)\right|}{\left|\delta\right|^{\gamma}}.
\end{equation}
 $\Lambda_{\gamma}(\bR^{n})$, the inhomogeneous Lipschitz space of order
 $\gamma$, is defined
as the space of continuous functions $f$ such that $
\left\|f\right\|_{\Lambda_{\gamma}} <\infty$. It is well known
that $f\in\Lambda_{\gamma}(\bR^{n})$ is of class $\mathcal
C^{k}(\bR^{n})$, with  $k:=\lfloor \gamma\rfloor$. Moreover,  for
$\alpha$ a multi-index with $|\alpha|\leq k$,
\begin{equation}\label{deriv}
\left\|\partial^{\alpha}f\right\|_{\Lambda_{\gamma-|\alpha|}}\leq
C(n,\gamma)\left\|f\right\|_{\Lambda_{\gamma}}.
\end{equation}
Similarly, we define the homogeneous Lipschitz
$\dot{\Lambda}_{\gamma}(\bR^{n})$ with
\begin{equation}
\left\|f\right\|_{\dot{\Lambda}_{\gamma}}=\sup_{x\in\bR^{n}}
\sup_{\delta\in\bR^{n}\setminus\left\{0\right\}}\frac{\left|D^{k+1}_{\delta}f(x)\right|}{\left|\delta\right|^{\gamma}},
\end{equation}

\section{Proofs of Theorem \ref{main0} and Theorem \ref{main}}
\begin{proof}[Proof of Theorem \ref{main0}] To simplify notations, we will
write $\H ^p$ in place of $\H^p(\bR^n)$, $\BMO$ in place of
$\BMO(\bR^n)$, etc.. The proof is inspired by the one given in
\cite{BIJZ} for the product $b\times h$ when $b$ is in $\BMO$.
Recall that we assume that $b\in\bmo$. The function $h\in\H^1$
admits an atomic decomposition with bounded atoms,
$$h:=\sum_j\lambda_j a_j\ \ \ , \ \ \sum_j|\lambda_j|\leq C\|h\|_{\H^1}.$$
When the sequence $h_\ell$ tends to $h$ in $\H^1$, the product
$b\times h_\ell$ tends to $b\times h$ as a distribution. So  we
can write
$$b\times h=\sum_j \lambda_j (b\times a_j),$$
where the limit is taken in the distribution sense. Since the
$a_j$ are bounded functions with compact support, the product
$b\times a_j$ is given by the ordinary product. We want to write
$b\times h:=h^{(1)}+h^{(2)}$, with $h^{(1)}\in L^1$ and
$h^{(2)}\in\H^\Phi_*$. Let us write, for each term $a_j$, which is
assumed to be adapted to $B_j$,
$$b\times a_j= (b-b_{B_j})a_j+b_{B_j}a_j.$$
By the $\BMO$ property as well as the fact that $|a_j|\leq
|B_j|^{-1}$, we have the inequality
$$\sum_j|\lambda_j|\int_{\bR^n}|b-b_{B_j}|\,|a_j|dx\leq C
\|b\|_{\bmo}\,\|h\|_{\H^1}.$$ Here $\left\|b\right\|_{\bmo}$ is
the sum of the two finite quantities that appear in the definition
of $\bmo$ given by \eqref{bmo}. We call
$$h^{(1)}:=\sum_j \lambda_j(b-b_{B_j})a_j,$$
which is the sum of a normally convergent series in $L^1$. Since
convergence in $L^1$ implies  convergence in the distribution
sense, it follows that $h^{(2)}$ is
$$h^{(2)}:= \sum_j \lambda_j b_{B_j} a_j, $$
 which is well defined in the distribution sense. Moreover
\begin{eqnarray*}
\M h^{(2)}&\leq &\sum_j |\lambda_j| |b_{B_j}|\M a_j\\
& \leq &\sum_j |\lambda_j| |b-b_{B_j}|\M a_j+ |b|\sum_j
|\lambda_j| \M a_j .\end{eqnarray*} The first term is in $L^1$
since $\M a_j\leq |B_j|^{-1}$. In order to conclude, we have to
prove that the second term is in $L^\Phi_*$. We first use the fact
that $\|\M a_j\|_1\leq C$ for some uniform constant $C$, which is
classical and may be found in \cite{St} for instance. Then we have
to prove that, for $\psi\in L^1$, the product $b\psi$ is in
$L^\Phi_*$. We claim that $b$ belongs uniformly to the exponential
class on each ball of measure $1$. Indeed, by John-Nirenberg
Inequality which is valid for $b$, for some constant $C$, which
depends only on the dimension, and for each ball $B$ such that
$|B|=1$,
\begin{equation}\label{JN}
\int_B\;\exp\left(\frac{\left| b(x)- b_B\right|}{C \|
b\|_{\bmo}}\right)\; dx\;\leqslant 2\;.
\end{equation}
Moreover, since $b\in\bmo$,  we have the inequality $|b_B|\leq
\|b\|_{\bmo}$. To prove that $b\psi$ is in $L^1$, we first
consider each such ball separately.   We use the following lemma,
which is an adaptation of lemmas given in \cite{BIJZ}.
changed
\begin{lem} If the integral on $B$ of $\exp |b|$  is bounded by $2$,
then, for some constant $C$,
$$\|b\psi|\|_{L^\Phi(B)}\leq C  \int_B |\psi|dx.$$
\end{lem}
\begin{proof}
 By homogeneity it is sufficient to find some constant $c$ such  that, for $\int_B |\psi|dx=c$  we have
 $$\int_B \frac {|b\psi|}{\log (e+|b\psi|)}dx\leq 1.$$ If we cut the integral into two parts depending on the fact
that $|b|<1$ or not, we conclude directly that the first part is bounded by $c$,
since we have a majorant by suppressing the denominator. For the
second part, we can suppress $b$ in the denominator. Then, we use
the duality between the $L \log L$ class, and the Exponential
class. It is sufficient to prove that the Luxembourg norm of
$\frac {|\psi|}{\log (e+|\psi|)}$ in the class $L \log L$  is
bounded by $1/2$ for $c$ small enough, which is elementary.
\end{proof}

We have an estimate for each cube $j+\mathbb Q$, and sum up.
This finishes the proof of Theorem \ref{main0}.


\end{proof}

Since $\bmo$ is the dual of $\h^1$, it is natural to see what is
valid for $h\in\h^1$. We can state the following.
\begin{thm}\label{local0}
For $h$ a function in $\h^1(\bR^n)$ and $b$ a function in $\bmo
(\bR^n)$, the product $b\times h$ can be given a meaning in the
sense of distributions. Moreover, we have the inclusion
\begin{equation}\label{inclusion1}
  b\times h\in L^1(\bR^n)+
  \h^\Phi_*(\bR^n).
\end{equation}
\end{thm}
\begin{proof} Again, we start from the atomic decomposition of $h$. In view of \eqref{add}, it is sufficient
 to  consider only those atoms $a_j$ that are adapted to balls $B$ such that
 $|B|\geq 1$. Remember that  they do not satisfy the moment
 condition \eqref{moment}. This one was only used to insure that
 $\|\M a_j\|_1\leq C$ for some independent constant. We now have $\|\M^{(1)}a_j\|_1\leq
 C$ since $\M^{(1)}a_j$, which is bounded by $|B_j|^{-1}$ is supported in the ball of same center as
 $B_j$ and radius twice the radius of $B_j$. Except for this
 point, the proof is identical.
\end{proof}

 Before leaving the case $p=1$, let us add some
remarks. Multipliers of the space $\BMO$ have been characterized
by Stegenga in \cite{S} (see also \cite{CL}) when $\bR^n$ is
replaced by the torus. It is easy to extend this characterization
to $\bmo$. Let us first define the space $\lmo$ as the space of
locally integrable functions $b$ such that
\begin{equation}\label{lmo}
\sup_{|B|\leq 1}\left( \frac {\log(e+1/|B|)} {|B|}\int_B
|b-b_B|dx\right)<\infty\ \ \ \ \mbox{\rm and } \sup_{|B|\geq
1}\left( \frac 1 {|B|}\int_B |b|dx\right)<\infty.
\end{equation}
$\lmo$ stands for {\sl logarithmic mean oscillation}.
\begin{prop}\label{mult-bmo} Let $\phi$ a locally integrable
function. Then the following properties are equivalent.
\begin{enumerate}
    \item [\rm{(i)}] The function $\phi$ is bounded and belongs to the space
    $\lmo$.
    \item[\rm{(ii)}] For every $b\in\bmo$, the function $b\phi$ is
    in $\bmo$.
\end{enumerate}
\end{prop}
\begin{proof}
We give a direct proof, which is standard,  for completeness. The
proof of (i)$\Rightarrow $(ii) is straightforward. Indeed, let us
first consider balls $B$ such that $|B|\leq 1$.  Writing $b=
(b-b_B) +b_B$, we conclude directly for the first term, and have
to prove that
$$|b_B |\times \frac {1} {|B|}\int_B
|\phi-\phi_B|dx \leq C \|b\|_{\bmo}\|\phi\|_{\lmo}.$$ Let $B'$ the
ball of same center as $B$ and radius $1$. It is well known that
the fact that $b$ is in $\BMO$ implies that $$|b_B-b_{B'}|\leq C
\log(e+1/|B|) \|b\|_{\BMO}.$$ We conclude, using the fact that
$|b_{B'}|\leq C\|b\|_{\bmo}. $ The proof  is even simpler for
balls $B$ such that $|B|\geq 1$.

Conversely, assume that we have (ii). Taking $b=1$, we already
know that $\phi$ is in $\bmo$. Also, by the closed graph theorem,
we know that there exists some constant $C$ such that, for every
$b\in\bmo$, the function
$$\|b\phi\|_{\bmo}\leq C\|b\|_{\bmo}.$$
We first claim that $\phi$ is bounded. By the Lebesgue
differentiation theorem, it is sufficient to prove that, for each
ball $B$, the mean $\phi_B$ is bounded. But $\phi_B=\la \phi,
|B|^{-1}\chi_B\ra =\la b\phi, a\ra$, where $a$ is some atom of
$\h^1$ and $b$ is bounded by $1$. Indeed, the characteristic
function $\chi_B$ may be written as the square of a function of
mean zero, taking values $\pm 1$ on $B$. So $\phi_B$ is bounded.
Now, since $\phi$ is bounded, the assumption implies that, for a
ball $B$ such that $|B|\leq 1$,
$$|b_B |\times \frac {1} {|B|}\int_B
|\phi-\phi_B|dx \leq C \|b\|_{\bmo}.$$ It is sufficient to find a
function $b$ with norm bounded independently of $B$ and such that
$|b_B |\geq c\log(e+1/|B|)$ The function $\log (|x-x_B|^{-1})$,
with $x_B$ the center of $B$, has this property.
\end{proof}
The previous proposition allows an interpretation in view of Theorem
\ref{main0}. The duals of Hardy-Orlicz spaces have been studied by
S. Jansen \cite{J}, see also the work of Viviani \cite{V} where
duality is deduced from their atomic decomposition. In particular,
the dual of the space $\h^{\Phi}$ is the space $\lmo$. It follows
that the dual of the space $L^1+ \h^{\Phi}$ is the space
$L^\infty\cap \lmo$. So if a duality argument was possible, which is not the case since we are not dealing with Banach
spaces, we would conclude that multiplication by $\bmo$ maps $\h^1$ into $L^1+ \h^{\Phi}$. Recall that we have a weaker statement.
\medskip

\begin{proof}[Proof of Theorem \ref{main}] When $p>\frac n{n+1}$,
the proof is an easy adaptation of the previous one. We start
again from an atomic decomposition of $h$ and define $h^{(1)}$ and
$h^{(2)}$ as before. To conclude for $h^{(1)}\in L^1$, it is
sufficient to prove that, for all balls $B$, one has
$$\int_B |b-b_B|dx\leq C|B|^{\frac 1p}\|b\|_{\Lambda_\gamma}.$$
If   $B$ has center $x_B$ and radius $r$, it follows at once from
the inequality $|b(x)-b(x_B)|\leq r^\gamma \|b\|_{\Lambda_\gamma}
\leq |B|^{\gamma/n}\|b\|_{\Lambda_\gamma}$, and the choice
$\gamma=n(1/p-1)$.

Next we conclude directly for $h^{(2)}$, using the fact that $b$
is bounded, so that $\M h^{(2)}\leq \|b\|_\infty \sum_j
|\lambda_j| \M a_j$. This last quantity is in $L^p$ since  $$\int
|\M h^{(2)}|^p \leq \|b\|_\infty \sum_j |\lambda_j|^p \int |\M
a_j|^p$$ and $\M a_j$'s are uniformly in $L^p$.

\medskip
For smaller values of $p$, we start again from an atomic
decomposition of $h$, but choose the atoms $a_j$ to be $(p,
\infty, s)$ for $s$ to be chosen later, that is, to satisfy the
moment condition \eqref{moment} up to order $s$. We then have to
modify the choice of $h^{(1)}$ and $h^{(2)}$ in order to be able
to treat the first term as above. We use the following definition.
\begin{defn} For $f$ a locally square integrable function
and $B$ a ball in $\bR^n$, we define $P^{k}_{B}f$ as the
orthogonal projection in $L^2(B)$ of $f$ onto the space of
polynomials of degree  $\leq k$.
\end{defn}
The next lemma is classical. It is the easy part of the
identification  of Lipschitz spaces with spaces of
Morrey-Campanato, see \cite{C}. We give its proof for
completeness.
\begin{lem} Let $\gamma>0$ and $k\geq \gamma $. There exists a constant $C$ such that, for $f$ a function in
$\lip$ and  $B$  a ball in $\bR^{n}$, then
\begin{equation}\label{morrey}
\frac{1}{\left|B\right|}\int_{B}\left|f(x)-P^{k}_{B}f(x)\right|dx\leq
C \left\|f\right\|_{\Lambda_{\gamma}}
\left|B\right|^{\frac{\gamma}{n}}.
\end{equation}
\end{lem}
\begin{proof}
In fact we prove an $L^2$ inequality instead of an $L^1$, which is
better. In this case, it is sufficient to prove the same
inequality with $P^{k}_{B}f$ replaced by some polynomial $P$ of
degree $\leq k$. This allows to conclude for $\gamma$ not an even
integer. Indeed, take for $P$ the Taylor polynomial at point $x_B$
(assuming that $B$ has center $x_B$ and radius $r$) and order
$\lfloor \gamma \rfloor $, using the fact that it makes sense by
\eqref{deriv}. Then, by Taylor's formula, $|f-P|$ is bounded on
$B$ by $Cr^{\gamma}\leq C |B|^{\gamma/n}$. For $\gamma$ an
integer, we conclude for \eqref{morrey} by interpolation.
\end{proof}
Let us come back to the proof of Theorem \ref{main}. We start
again from an atomic decomposition of $f$,
 We fix $ k\geq \gamma$
and pose
$$h^{(1)}:=\sum_j \lambda_j(b-P^{k}_{B_j}b)a_j.$$
Using the previous lemma, we conclude as before that $h^{(1)}$ is
in $L^1$. In order to have $h^{(2)}:=b\times h-h^{(1)}$  in
$\H^p$, it is sufficient that each term $(P^{k}_{B_j}b) a_j$ be,
up to the multiplication by a uniform constant, a $(p, \infty,
s')$-atom with $s'\geq\gamma$. The moment condition is clearly
satisfied if we have $s\geq k+s'$. We can in particular choose
$k=s'=\lfloor\gamma\rfloor$ and $s=2\lfloor\gamma\rfloor$. It
remains to prove that $P^{k}_{B_j}$ is uniformly bounded. This
follows from the following lemma.
\begin{lem}\label{bornitude de la projection}
Let $k$ be a positive integer. There exists a constant $C>0$ such
that for every ball $B$ in $\bR^{n}$,
\begin{equation}
\left\|P^{k}_{B}f\right\|_{L^{\infty}(B)}\leq C
\left\|f\right\|_{L^{\infty}(B)},
\end{equation}\label{norme proj de f}
for all functions $f$ which are bounded on the ball $B$.
\end{lem}
\begin{proof}
We remark first that, by invariance by translation we can assume
that $B$ is centered at $0$. Next, by invariance by dilation, we
can also assume that $|B|=1$. So we have to prove it for just one
fixed ball. Now, since the  projection is done on a  finite
dimensional space,
$$
\left\|P^{k}_{B}f\right\|_{L^{\infty}(B)}\leq C_k
\left\|P^{k}_{B}f\right\|_{L^{2}(B)} \leq
C_k\left\|f\right\|_{L^{2}(B)}\leq C_k
\left\|f\right\|_{L^{\infty}(B)}.
$$
\end{proof}
This allows to conclude for the proof of the theorem.
\end{proof}
As for the case $p=1$, we can take $h$ in the local Hardy space.
\begin{thm}\label{local}
For $h$ a function in $\h^p(\bR^n)$ and $b$ a function in
$\lip$,
 the product $b\times h$ can be given a meaning in the
sense of distributions. Moreover, we have the inclusion
\begin{equation}\label{inclusion2}
  b\times h\in L^1(\bR^n)+ \h^p(\bR^n).
\end{equation}
\end{thm}
\begin{proof}  The adaptation of
the previous proof is done in the same way as we have done for
Theorem \ref{local0} compared to Theorem \ref{main0}. We leave it
 to the reader.
\end{proof}
We did not give estimates of the norms, but it follows from the
proof of Theorem \ref{main} that we have the inequality
$$\|h^{(1)}\|_1 + \|h^{(2)}\|_{\H^p}\leq C \|h\|_{\H^p}\times
\|b\|_{\lip}.$$ So the bilinear operator
\begin{eqnarray*}
\mathfrak P:& \lip\times\H^{p}(\bR^{n})\rightarrow  &L^{1}(\bR^{n})+\H^{p}(\bR^{n})\\
&\left(b,h\right) \mapsto  & b\times h
\end{eqnarray*}
is continuous. It is easy to see that the term in $L^1$ is present
in general: for instance take an example in which the product is
positive. The same remarks are valid for all three other cases.
\begin{remark} The product of $b\in
\dot{\Lambda}_{\gamma}(\bR^{n})$ and $h\in \H^{p}(\bR^{n})$ is
also well defined. It belongs to some
$L^{1}(\bR^{n})+\H^{p}_w(\bR^{n})$ with a weight $w$ conveniently
chosen. Now $b$ is no more bounded, but can increase as
$|x|^\gamma$ at infinity. We can take any weight
$(1+|x|)^{-\alpha}$, with $\alpha>\gamma p$.
\end{remark}

\section{Generalization to spaces of homogeneous type}

All proofs generalize easily to spaces of homogeneous type once
one has been able to define correctly the product $b\times h$. We
will not give into details of terminology and proofs when the
generalization may be done without any  difficulty, but will
essentially concentrate on the definition of the product.

Let us first recall some definitions. We assume that we are given
a locally compact Hausdorff space $X$, endowed with a quasi-metric
$d$ and a positive regular measure $\mu$ such that the doubling
condition
\begin{equation}
0<\mu \left( B_{\left(x,2r\right) }\right) \leq C\mu
\left(B_{\left( x,r\right)}\right) <+\infty  \label{doublement}
\end{equation}
holds, for all $x$ in $X$ and $r>0.$ Here, by a quasi-metric $d$,
we mean  a function $d:X\times X\rightarrow \left[ 0;+\infty
\right[ $ which satisfies
\begin{enumerate}
\item[( i)] $d(x,y)=0$ if and only if $x=y$ ;

\item[( ii)] $d(x,y)=d(y,x)$ for all $x,y$ in $X$;

\item[( iii)] there exists a finite constant $\kappa \geq 1$
such that
\begin{equation}
d(x,y)\leq \kappa \left( d(x,z)+d(z,y)\right)  \label{quasimetric}
\end{equation}
for all $x,y,z$ in $X.$
\end{enumerate}
Given $x\in X$ \ and $r>0,$ we note $B_{\left( x,r\right)
}=\left\{ y\in X:d\left( x,y\right) <r\right\} $  the ball with
center $x$ and radius $r.$
\begin{defn}
We call space of homogeneous type $\left( X,d,\mu \right) $ such a
locally compact space $X$, given together with the quasi-metric
$d$ and the nonnegative Borel measure $\mu $ on $X$ that satisfies
the doubling condition.
\end{defn}

On such a space of homogeneous type we can define the $\BMO(X)$
and $\bmo(X)$ spaces and have the John-Nirenberg inequality. Just
replace Euclidean balls by the balls on $X$, and the Lebesgue
measure $dx$ by the measure $d\mu$. The Hardy space $\H^1(X)$ can
be defined by the atomic decomposition and we have the duality
$\H^1$-$\BMO$. But we want also  to define the product $b\times h$
for $b\in \bmo(X)$ and $h\in \H^1(X)$, while we can no more speak
of distributions. In view of Theorem \ref{mult-bmo}, which
generalizes easily in this context, we can take $\C\cap
L^{\infty}\cap \lmo$ as a space of test functions, but we need to
have density theorems of such functions in the space of continuous
compactly supported functions to recover the pointwise product
when it is integrable. We encounter a fortiori this difficulty
when dealing with Lipschitz spaces.\medskip

Mac\'ias and Segovia have overcome this kind of difficulty for
being able to develop the theory of $\H^p$ spaces for $p<1$ (see
also \cite{U}). We will assume that the measure $\mu$ does not
charge points for simplification. They have proved \cite{ms1}
that, without loss of generality, $X$ may be assumed a {\sl normal
space of homogeneous type and of order} $\alpha>0$ when,
eventually, the quasi-distance is replaced by an equivalent one.
That is, the quasi-metric $d$ and the measure $\mu $ are assumed
to satisfy the following properties.

\noindent There exist four positive constants $A_{1},A_{2},K_{1}$
and $K_{2}$, such that
\begin{eqnarray}
    A_{1}r\leq\mu(B_{(x,r)})\leq A_{2}r &if& 0\leq r\leq K_{1}\mu (X)\\
    B_{(x,r)}=X &if& r>K_{1}\mu(X).
\end{eqnarray}
\begin{equation}
    \left|d(x,z)-d(y,z)\right|\leq K_{3}r^{1-\alpha}d(x,y)^{\alpha},\label{normal}
\end{equation}
for every $x,y$ and $z$ in $X$, whenever $d(x,z)<r$ and $d(y,z)<r$.

\smallskip

In  the Euclidean case, a  normal quasi-distance   is given by
$|x-y|^n$. The assumption \eqref{normal} is satisfied with
$\alpha=1/n$.

Let us then define the Lipschitz spaces, as it is natural, by the
following.
\begin{defn}
Let $\gamma>0$. The Lipschitz space $\Lambda_{\gamma}(X,d,\mu)$
consists of bounded continuous  functions $f$ on $X$ for which,
for some constant $C$ and all $x,y$,
\begin{equation}
\left| f(x)-f(y)\right|\leq C d(x,y)^{\gamma}. \label{lips}
\end{equation}
\end{defn}
Note that there is a change of parameter in the Euclidean space,
$\gamma$ has been changed into $\gamma/n$.

Remark that with these definitions Lipschitz spaces are contained
in $\lmo$.

 We know that the space
is not reduced to $0$ when $\gamma$ is not larger than $\alpha$
because of the fact that the distance itself satisfies this kind
of condition. Mac\'ias and Segovia have proved in \cite{ms2} that
one can build approximate identities in order to approach
continuous functions with compact support by Lipschitz functions
of order $\gamma$, for any $\gamma<\alpha$.   Moreover, they
define the space of distributions $\left(E^{\alpha}\right)^{\ast}$
as the dual of the space $E^{\alpha}$, consisting of all functions
with bounded support, belonging to $\Lambda_{\beta}$ for every
$0<\beta<\alpha$. From this point, they can use distributions to
define $\H^p$ spaces when $p>(1+\alpha)^{-1}$. We recover in the
Euclidean case the condition $p>n/(n+1)$, which is the range where
atoms are assumed to satisfy only the moment condition of order
zero, and where the dual is a Lipschitz space defined by a
condition implying only one difference operator.

This notion of distribution is exactly what we need for the
definition of products. With the conditions above and in the
corresponding range of $p$, products may be defined in the
distribution sense and the four theorems are valid.

\medskip

Remark that, for $X$ the boundary of a bounded smooth
pseudo-convex domain of finite type, Lipschitz spaces can be
defined for all values of $\gamma$ and $\H^p$ spaces can  be
defined for $p$ arbitrarily small. The moment conditions of higher
order rely on the use of vector fields related to the geometric
structure of the boundary. We refer to \cite{McN} for the
geometrical aspects, and to \cite{BG} for the detailed statements
related to  products of holomorphic functions in $\H^p$ and $\BMO$
or Lipschitz spaces.


\begin{thebibliography}{MTW1}

\bibitem[BG]{BG} A. B{\scshape onami} and S. G{\scshape rellier},
{\em Decomposition theorems for Hardy-Orlicz spaces and weak
factorization}, preprint.

\bibitem[BIJZ]{BIJZ} A. B{\scshape onami}, T. I{\scshape waniec}, P. J{\scshape ones}, M. Z{\scshape insmeister},
 { \em On the product of functions in $BMO$ and $\mathcal H^1$},  Ann. Inst. Fourier, Grenoble {\bf 57}, 5 (2007) 1405-1439.


\bibitem [C]{C}  S. C{\scshape ampanato}, {\em Propriet\`a di h\"olderianit\`a di alcune classi di funzioni},
Ann. Scuola Norm. Sup. Cl. Sci. {\bf 17} (1963), 175-188.

\bibitem[CL]{CL} D. C. C{\scshape hang} and S. Y. L{\scshape i}, {\em On the boundedness of multipliers,
commutators and the second derivatives of Green's operators on
$H^{1}$ and $BMO$}, Ann. Sc. Norm. Super. Pisa, Cl. Sci., IV. Ser.
{\bf 28} $n^{\circ}2$ (1999), 341-356.

\bibitem[CW1]{CW1} R. C{\scshape oifman}, G. W{\scshape eiss},
{ \em Analyse Harmonique Non-commutative sur Certains Espaces
Homog`enes},   (1977), 569-645. Lecture Notes in Math. {\bf 242},
Springer, Berlin, 1971.

\bibitem[CW2]{CW2} R. C{\scshape oifman}, G. W{\scshape eiss},
{ \em Extensions of Hardy spaces and their use in analysis},
Bull. Amer. J. Math. {\bf 83} (1977), 569-645.


\bibitem[G]{G} D. G{\scshape olberg}, {\em A local version of Hardy spaces},
{ Duke J. Math.} {\bf 46} (1979), 27-42.


\bibitem [GLY]{GLY} L. G{\scshape rafakos}, L. L{\scshape iu} and D. Y{\scshape ang}, {\em Maximal function characterizations of Hardy spaces on RD-spaces and their
applications}, preprint.

\bibitem[Gr]{g} L. G{\scshape rafakos}, {\em Classical and Modern Fourier Analysis}, Pearson Education, Inc. Upper Saddle River, New Jersey 07458.

 \bibitem[J]{J} S. J{\scshape anson}, {\em Generalizations of Lipschitz spaces and an application to Hardy spaces and bounded mean
 oscillation},
   Duke Math. J. {\bf  47}  (1980), no. 4, 959--982.

\bibitem[JPS]{JPS} S. J{\scshape anson}, J. P{\scshape eetre} and S. S{\scshape emmes},
{\em On the action of Hankel and Toeplitz operators on some function spaces}, Duke Math. J. {\bf 51} (1984),
937--958.

\bibitem[KL]{KL} S. G. K{\scshape rantz} and S.-Y. L{\scshape i},
{\em Boundedness and Compactness of Integral Operators on Spaces
of Homogeneous Type and Applications}, J. Math. Anal. Appl. {\bf
258} No.2 (2001), 642-657.


\bibitem[McN]{McN} J. D. McN{\scshape eal}, {\em Estimates on the Bergman kernels of
convex domains}, { Adv. in Math}  {\bf 109} (1994), 108-139.


\bibitem [MS1]{ms1}R. A. M{\scshape ac$\acute{\text{i}}$as} and C. S{\scshape egovia}, { \em Lipschitz function on Spaces of Homogeneous type}, Advances in Math. {\bf 33} (1979), 257-270.

\bibitem [MS2]{ms2}R. A. M{\scshape ac$\acute{\text{i}}$as} and C. S{\scshape egovia}, { \em A decomposition into Atoms of Distributions on Spaces of Homogeneous type}, Advances in Math. {\bf 33} (1979), 271-309.

\bibitem[S]{S} D. A. S{\scshape tegenga}, {\em Bounded Toeplitz operators on $H\sp{1}$ and
applications of the duality between $H\sp{1}$ and the functions of
bounded mean oscillation}, Amer. J. Math.  {\bf 98}  (1976), no.
3, 573--589.

\bibitem[St]{St} E. M. S{\scshape tein}, {\em Harmonic analysis, real-variable
methods, orthogonality, and oscillatory integrals},
Princeton Math. Series 43, Princeton
University
Press, Princeton 1993.

\bibitem [U]{U} A. U{\scshape chiyama} {\em The factorization of $H\sp{p}$ on the space of homogeneous type},
Pacific J. Math.  {\bf 92}  (1981),  453--468.
\bibitem[V]{V} B. E. V{\scshape iviani}, {\em An Atomic Decomposition of the Predual
 of BMO$(\rho)$}, Revista Matem$\acute{\text{a}}$tica Iberoamericana {\bf 3}  (1987), 401-425.

\end{thebibliography}
\end{document}